\newif\ifextended\extendedfalse

\documentclass{elsarticle_preprint}
\usepackage[utf8]{inputenc}
\usepackage[english]{babel}
\usepackage{amsmath,amssymb,amsfonts,amsthm}
\usepackage[disable]{todonotes}
\usepackage{comment}
\usepackage[binary-units]{siunitx}     
\usepackage{enumerate}
\usepackage{url}
\usepackage{bm}
\usepackage{pgf,pgfplots,pgfplotstable}
\usepackage[linesnumbered, ruled, noline, noend]{algorithm2e}
\usepackage{mathtools}
\usepackage{hyperref}
\usepackage{cleveref}
\usepackage{caption}
\usepackage{placeins}
\usepackage{subcaption}
\pgfplotsset{compat=1.16}

\usetikzlibrary{arrows.meta, 
                chains,
                positioning,
                shapes}    
\tikzset{reset join/.code={\def\tikz@after@path{}}}

\theoremstyle{definition}
\newtheorem*{remark}{Remark}
\theoremstyle{plain}
\newtheorem{theorem}{Theorem}
\newtheorem{lemma}{Lemma}

\newcommand{\sig}[2][]{\normalfont\text{Sig}\left(#1\bm{#2}\right)}
\newcommand{\gen}[1]{\normalfont\text{Gen}\left(#1\right)}
\newcommand{\flag}[1]{\normalfont\text{Flag}\left(#1\right)}
\newcommand{\tail}[1]{\normalfont\text{Tail}\!\left(#1\right)}
\newcommand{\algo}{\texttt{M5GB}}
\newcommand{\reducealgo}{\texttt{Reduce}}
\newcommand{\updatealgo}{\texttt{Update}}
\DeclareMathOperator{\regmod}{\textnormal{mod}_{reg}}

\begin{document}

\title{A Signature-Based Gröbner Basis Algorithm with Tail-Reduced Reductors (\texttt{M5GB})}
\date{}
\author[1]{Manuel Hauke}
\ead{hauke@math.tugraz.at}

\author[2]{Lukas Lamster}
\ead{lukas.lamster@iaik.tugraz.at}

\author[2]{Reinhard Lüftenegger}
\ead{reinhard.lueftenegger@iaik.tugraz.at}

\author[2]{Christian Rechberger}
\ead{christian.rechberger@iaik.tugraz.at}

\address[1]{Institute of Analysis and Number Theory, Graz University of Technology, Steyrergasse 30/II, 8010 Graz, Austria}
\address[2]{Institute of Applied Information Processing and Communications, Graz University of Technology, Inffeldgasse 16a, 8010 Graz, Austria}

\tnotetext[t1]{Author list in alphabetical order; see \url{https://www.ams.org/profession/leaders/CultureStatement04.pdf}.}

\begin{abstract}
Gröbner bases are an important tool in computational algebra and, especially in cryptography, often serve as a boilerplate for solving systems of polynomial equations. Research regarding (efficient) algorithms for computing Gröbner bases spans a large body of dedicated work that stretches over the last six decades. The pioneering work of Bruno Buchberger in 1965 can be considered as the blueprint for all subsequent Gröbner basis algorithms to date. Among the most efficient algorithms in this line of work are signature-based Gröbner basis algorithms, with the first of its kind published in the late 1990s by Jean-Charles Faugère under the name \texttt{F5}. In addition to signature-based approaches, Rusydi Makarim and Marc Stevens investigated a different direction to efficiently compute Gröbner bases, which they published in 2017 with their algorithm \texttt{M4GB}. The ideas behind \texttt{M4GB} and signature-based approaches are conceptually orthogonal to each other because each approach addresses a different source of inefficiency in Buchberger's initial algorithm by different means.

We amalgamate those orthogonal ideas and devise a new Gröbner basis algorithm, called \algo{}, that combines the concepts of both worlds. In that capacity, \algo{} merges strong signature-criteria to eliminate redundant S-pairs with concepts for fast polynomial reductions borrowed from \texttt{M4GB}. We provide proofs of termination and correctness and a proof-of-concept implementation in \verb|C++| by means of the Mathic library. The comparison with a state-of-the-art signature-based Gröbner basis algorithm (implemented via the same library) validates our expectations of an overall faster runtime for quadratic overdefined polynomial systems that have been used in comparisons before in the literature and are also part of cryptanalytic challenges.
\end{abstract}

\maketitle

\section{Introduction}

Gröbner bases are an essential tool in commutative algebra and algebraic geometry. Several important applications in these areas are \textit{(a)} the \textit{Ideal Equality Problem}, characterizing the equality of two ideals through the reduced Gröbner bases of their sets of generators, \textit{(b)} the \textit{Ideal Membership Problem}, characterizing whether a polynomial belongs to a given ideal via the division remainder modulo the respective (reduced) Gröbner basis, and \textit{(c)} the \textit{Elimination Problem}, eliminating variables from a system of polynomial equations through, e.g., lexicographic Gröbner bases. There are many more applications of Gröbner bases in signal and image processing, robotics, automated geometric theorem proving, and solving systems of polynomial equations \cite{BW98}, \cite{BS10}.

Especially in cryptography, both public-key and symmetric cryptography, the problem of solving systems of polynomial equations arises in different contexts, ranging from block cipher and hash function analysis \cite{ACGKLRS19,BPW06} to the analysis of asymmetric encryption and signature schemes \cite{TPD21,FJ03}.
The applicability of Gröbner bases to cryptanalysis has been one of the driving factors for research on efficient algorithms for computing them.
In general, the pioneering \texttt{Buchberger} algorithm for computing Gröbner bases devised by Bruno Buchberger and published in 1965 \cite{B65,B76} can be considered highly inefficient, mainly due to the excessive amount of redundant computations that do not provide any new information for the eventual Gröbner basis. In more detail, the \texttt{Buchberger} algorithm repeatedly reduces so-called \textit{S-pairs} \cite[p.84]{CLO02}, adds all non-zero remainders to the current basis and repeats this process until all S-pairs reduce to zero with respect to the current basis. During the process of repeatedly reducing S-pairs, often many of those S-pairs reduce to zero and thus they do not provide any new information. To tackle this inefficiency, further criteria have been developed to streamline the \texttt{Buchberger} algorithm by detecting and discarding S-pairs that would otherwise reduce to zero. The work by Gebauer and Möller \cite{GM88} implements these criteria and presents a more efficient instantiation of the \texttt{Buchberger} algorithm.

A different approach, initially investigated by Buchberger \cite{B83,B83a} and Lazard \cite{L79,L83,L01} and further developed by Faugère in 1999 \cite{F99}, relates the problem of reducing S-pairs to the problem of reducing matrices. The basic underlying idea is that for a degree bound large enough and all term multiples of the initial polynomials up to this degree, the matrix containing the corresponding coefficients of the term multiples yields, after Gaussian row reduction, a Gröbner basis of the ideal generated by the initial polynomials. Rather than choosing a degree bound large enough and constructing one large matrix, Faugère's \texttt{F4} algorithm in \cite{F99} constructs matrices for smaller degrees, row-reduces the corresponding smaller matrices and continues in this fashion until a Gröbner basis is found. Compared to the \texttt{Buchberger} algorithm, the advantage of \texttt{F4} is that S-pairs are reduced in parallel rather than sequentially. This advantage is the main source of the particular efficiency of \texttt{F4} and some of the fastest Gröbner basis implementations to date rely on this approach, as, e.g., the implementation in the computer algebra system Magma.

Signature-based Gröbner basis algorithms are another line of work regarding more efficient instantiations of the \texttt{Buchberger} algorithm. In the \texttt{F5} algorithm introduced by Faugère in 2002 \cite{F02}, so-called \textit{signatures} help to keep track from which initial polynomials some S-pair has been calculated. The information from the signatures allows to detect whether a S-pair reduces to zero without having to carry out the reduction. Thus, the main idea of signature-based criteria is reducing the amount of redundant reductions. For a particular class of polynomial systems, called regular sequences, \texttt{F5} does not carry out any redundant reduction at all. The \texttt{F5} algorithm and other signature-based Gröbner basis algorithms have later been incorporated into the \textit{rewrite} framework published by Christian Eder and Bjarke Roune \cite{ER13}. The rewrite framework generalizes many different (signature-based) approaches for computing Gröbner bases and unifies them under the umbrella of a single comprehensive framework.

Compared to the approaches discussed above, Rusydi Makarim and Marc Stevens present an orthogonal concept for computing Gröbner bases \cite{MS17}: their \texttt{M4GB} algorithm is based on the Gebauer-Möller version of the \texttt{Buchberger} algorithm, with the difference that reductions of S-pairs are carried out with tail-reduced reductors. In addition, these tail-reduced reductors are stored for potential later reuse. The main advantage of this approach is that no new reducible terms are introduced into the reduction remainder, which allows reducing S-pairs in a term-wise and recursive manner. This results in a fast polynomial reduction routine. The downside, however, is that any reductor has to be tail-reduced first, and the overall advantage of this approach depends on how often the algorithm is able to reuse already constructed and stored tail-reduced reductors.

\subsection{Our Contribution}

We present a new algorithm for computing Gröbner bases, called \algo{}, which combines the strengths of \texttt{M4GB} \cite{MS17} and signature-based Gröbner basis algorithms like \texttt{F5} \cite{F02}. We provide proofs of termination and correctness for \algo{}. In particular, we show how one can adapt the fast reduction routine used in \texttt{M4GB} to work with the signature-based criteria from \texttt{F5}-like algorithms. This creates a generic optimization that can be used for any signature-based Gröbner basis algorithm that does not use a matrix approach for polynomial reduction. The question of merging the fast reduction routine from \texttt{M4GB} with signature-based criteria arises naturally but resolving it is a non-trivial task that requires technical care, especially when it comes to algorithmic efficiency. To date, and to the best of our knowledge, no algorithm in this line of work has been published yet.

For a proof-of-concept implementation, we concentrate on the signature-based algorithm \texttt{SB} from Stillman and Roune \cite{RS12}, also called \texttt{SigGB} in the reference implementation \cite{mathicgb}, and adapt this algorithm to be compatible with an \texttt{M4GB}-like reduction routine. We show that using the same library for implementing \texttt{SB} and \algo{}, we obtain a significant, scalable speed-up for dense, quadratic, overdefined systems. These systems are used for benchmark purposes in the original article about \texttt{M4GB} by Makarim and Stevens \cite{MS17} and are posed as a problem instance in the MQ Challenge \cite{mq}.

\subsection{Related Work}
Compared to our approach for computing Gröbner bases in \algo{}, there exist related but different approaches in the literature. Here, we briefly discuss the main conceptual differences. One difference applies to all discussed algorithms below: we state our \algo{} algorithm in the rewrite framework \cite{ER13}, while the below algorithms adopt the basic structure of \texttt{F5}. The rewrite framework comprises the original \texttt{F5} algorithm as a special instantiation.

\paragraph{\texttt{F4\slash 5}} Albrecht and Perry \cite{AP10} describe an algorithm that combines \texttt{F4}-style reduction with \texttt{F5}-like signature criteria. This means, \cite{AP10} integrates - like \algo{} - a fast reduction routine with signature-based criteria to discard S-pairs. The difference to \algo{} is that their algorithm \texttt{F4\slash 5} uses the same linear algebra approach for reducing S-pairs as \texttt{F4}, and thus conceptionally resembles \texttt{Matrix-F5} \cite{BFS15} rather than \algo{}. Hence, it is the reduction of S-pairs that distinguishes \texttt{F4\slash 5} and \algo{}: the former algorithm uses \texttt{F4}-style reduction, while the latter one uses \texttt{M4GB}-style reduction. For a more detailed differentiation between the respective reduction routines in \texttt{F4} and \texttt{M4GB} we refer to \cite[Chapter 4.2]{MS17}. 

\paragraph{\texttt{F5C} and \texttt{F5R}} Eder and Perry \cite{EP10} present a variant of Faugère's original \texttt{F5} algorithm which works with reduced intermediate Gröbner bases rather than non-reduced ones. This results in fewer S-pairs to consider for checking the signature criteria and, eventually, fewer polynomial reductions. Eder and Perry differentiate their approach called \texttt{F5C}, ``F5 Computing by reduced Gröbner bases'', from the approach devised by Stegers \cite{S06} for which they use the denomination \texttt{F5R}, ``F5 Reducing by reduced Gröbner bases''. In \cite{S06}, Stegers' \texttt{F5R} algorithm uses reduced intermediate bases only for polynomial reductions, however, it still uses unreduced intermediate bases for computing new S-pairs. In contrast, \texttt{F5C} uses reduced intermediate bases for, both, polynomial reductions and generating new S-pairs. To summarize, the advantage of \texttt{F5R} over \texttt{F5} is faster polynomial reductions, while the advantage of \texttt{F5C} over \texttt{F5R} is a lower number of S-pairs to compute. The main conceptual difference between \texttt{F5C} and \algo{} is, again, the reduction routine: \texttt{F5C} uses ordinary polynomial reduction while \algo{} uses \texttt{M4GB}-style reduction.

\section{Preliminaries}

Any work treating the theory behind Gröbner bases and, in particular, describing different algorithms to compute Gröbner bases is faced with the challenge of having to introduce a significant amount of definitions and denominations. On top of that, there are often considerable notational differences between different authors. This being said, in \Cref{Sec:preliminary-definitions} we pay attention to stay close to commonly shared denominations and to find a balance between a rigorous and compact nomenclature.

Furthermore, in \Cref{Sec:m4gb} and \Cref{Sec:signature-algorithms} we give brief accounts of \texttt{M4GB} and signature-based algorithms for computing Gröbner bases, respectively. We describe these algorithms only to the extent that we are able to sketch their core ideas needed for our presentation of \algo{} in \Cref{Sec:new-algorithm}. We assume some familiarity with these algorithms from the reader, although the core ideas should become apparent without any deeper prior knowledge.

\subsection{Preliminary Definitions}
\label{Sec:preliminary-definitions}

We work with polynomials in the variables $X_1,\ldots,X_n$ over a finite field $\mathbb{F}$, i.e., with the polynomial ring $P\coloneqq\mathbb{F}[X_1,\ldots,X_n]$. A term is a power product of variables,
while a monomial is a product of coefficient and term. By $T$ we denote the set of all terms in $P$.
For a polynomial $f\in P$, the set $T(f)$ shall denote the set of all terms of $f$. For a term $t\in T(f)$, the corresponding coefficient is denoted as $C_t(f)$.
We define the free $P$-module $P^m$ with generators $\bm{e_1}\coloneqq (1,0,\ldots,0)$, \ldots, $\bm{e_m}\coloneqq (0,\ldots,0,1)$. As for polynomials, a \textit{module term} is an element in $P^m$ of the form $t\bm{e_i}$, while a \textit{module monomial} is an element of the form $c\cdot t\bm{e_i}$, for $c\in\mathbb{F}$, $t\in T$ and $1\leq i\leq m$. The set of all module terms in $P^m$ is denoted by $T_m$.

Throughout this article, we write module elements $\bm{f},\bm{g},\ldots$ in $P^m$ in boldface, whereas polynomials $f,g,\ldots$ in $P$ are written in normal style. We denote a term order on $T$ and a compatible order extension\footnote{An order extension is called compatible, if $\forall u,v\in T\;\forall 1\leq i\leq m: u\leq v\Rightarrow u\bm{e_i}\leq v\bm{e_i}$.} to module terms in $T_m$ by the same sign $\leq$. We believe, this ambiguity is justified by an easier notation and causes no deeper confusion because the context clarifies whether $\leq$ relates polynomials or module elements. For a given term order $\leq$, the \textit{leading term} of a polynomial $f\in P$, denoted by $LT(f)$, is defined as the $\leq$-maximum term in $T(f)$ and the \textit{leading coefficient} as the associated coefficient of $\text{LT}(f)$. In a similar fashion, the \textit{module leading term} $\text{MLT}(\bm{f})$ and \textit{module leading monomial} $\text{MLM}(\bm{f})$ are defined for a module element $\bm{f}\in P^m$ and a compatible order extension $\leq$. The polynomial $\text{Tail}(f)\coloneqq f-\text{LC}(f)\cdot\text{LT}(f)$ is called the \textit{tail of f}.

Given a finite set of non-zero polynomials $F\coloneqq\{f_1,\ldots,f_m\}\subseteq P\setminus\{0\}$, the module homomorphism $\varphi_F:P^m\rightarrow P$ given by $(p_1,\ldots,p_m)\mapsto \sum_{i} p_i f_i$ connects the module and polynomial perspective. Usually, the underlying set $F$ is clear, therefore we often omit the subscript and just write $\varphi$ instead of $\varphi_F$. Using the canonical generators of $P^m$, we can also write $\varphi:\sum_{i} p_i\bm{e_i}\mapsto \sum_{i} p_i f_i$. Any module element $\bm{h}\in P^m$ with $\varphi(\bm{h})=0$ is called a \textit{syzygy}. The \textit{signature} of a module element $\bm{f}\in P^m$ is given by $\sig{f}\coloneqq \text{MLT}(\bm{f})\in T_m$; of course, always relative to some compatible order extension $\leq$.

For a finite set of polynomials $G\subseteq P\setminus\{0\}$, a non-zero polynomial $f$ is said to be \textit{reducible with respect to $G$}, if there exist a term $t\in T(f)$ and an element $g\in G$ such that $\text{LT}(g) \mid t$. If $u\coloneqq t\slash\text{LT}(g)$ and $c\coloneqq C_t(f)$, we denote the reduction itself by
\(
f \longrightarrow_{G} f-c\cdot ug.
\)

The element $c\cdot ug$ is called a \textit{reductor of $f$}. If $t=\text{LT}(f)$, the reduction step is also called a \textit{top}-reduction, otherwise a \textit{tail-reduction} and the corresponding reductors are called \textit{top-reductor} and \textit{tail-reductor}, respectively. If a polynomial is not reducible (or tail-reducible) with respect to $G$, it is called \textit{irreducible} (or \textit{tail-irreducible}) \textit{with respect to $G$}. For the sake of notational convenience, any non-zero scalar multiple $d\cdot ug$, $d\in\mathbb{F}\setminus\{0\}$, is also called a reductor of $f$. This is why we often drop the scalar coefficient and just call $ug$ a reductor of $f$. If $f$ reduces to $h\in P$ in finitely many reduction steps with respect to $G$, we denote this by
\(
f\longrightarrow_{G,\ast}h.
\)
This also includes the case in which no reduction steps are done at all, hence
\(
f\longrightarrow_{G} f
\)
is trivially valid. We call a polynomial $f' \in P$ to be a \textit{normal form of $f$ with respect to $G$} if $f\longrightarrow_{G,\ast}f'$ and $f'$ is irreducible with respect to $G$. We use the denomination
\[
f\bmod G \coloneqq \{f' \in P: f' \text{ a normal form of } f \text{ w.r.t. } G\}
\]
to write down the set of all normal forms of a polynomial $f$.

\begin{remark}
Usually, we omit the specification \textit{with respect to $G$} and presume it to be clear from the context; whenever necessary, we explicitly mention the underlying set $G$. The same applies for $Sig$-reductions defined below. Furthermore, we often do not mention nor incorporate the underlying (module) term order in our definitions and terminology. Again, the aim is having a lighter notation.
\end{remark}

For a finite set of module elements $\bm{G} \subseteq P^m\setminus\{\bm{0}\}$, a non-zero module element $\bm{f}\in P^m\setminus\{\bm{0}\}$ is said to be \textit{$Sig$-reducible with respect to $\bm{G}$} if there exist a term $t \in T(\varphi(\bm{f}))$ and an element $\bm{g} \in \bm{G}$ such that the following two properties hold:
	\textit{(i)} $\text{LT}(\varphi(\bm{g})) \mid t$, in which case we set $u \coloneqq  t\slash LT(\varphi(\bm{g}))$;
	\textit{(ii)} $\sig{f} \geq \sig[u]{g}$.
If these properties are fulfilled, we define $\bm{f}-c\cdot u\bm{g}$ as the outcome of the $Sig$-reduction, where $c\coloneqq \text{C}_t(\varphi(\bm{f}))\slash\text{LC}(\varphi(\bm{g}))$, and denote the $Sig$-reduction itself by
\(
\bm{f} \longrightarrow_{\bm{G}} \bm{f}-c\cdot u\bm{g}.
\)
In particular, the element $c\cdot u\bm{g}$ is called a \textit{$Sig$-reductor of $\bm{f}$}. If $\sig{f} > \sig{u\bm{g}}$, we call it a \textit{regular $Sig$-reduction}, otherwise a \textit{singular $Sig$-reduction}. We denote a regular $Sig$-reduction by
\(
\bm{f} \longrightarrow_{\bm{G},reg} \bm{h} 
\)
and, analogously, any finite number of regular $Sig$-reductions on $\bm{f}$ to a module element $\bm{h}$ by
\(
\bm{f} \longrightarrow_{\bm{G},reg,\ast}\bm{h}.
\)
We say
\(
\bm{f} \longrightarrow_{\bm{G},\ast} 0
\)
if
\(
\bm{f} \longrightarrow_{\bm{G},\ast} \bm{h}
\)
for a syzygy $\bm{h}$. This notation is justified by $\varphi(\bm{h})=0$. We believe, the definitions of \textit{$Sig$-top-reduction}, \textit{$Sig$-tail-reduction}, \textit{$Sig$-irreducible}, \textit{$Sig$-tail-irreducible}, \textit{regularly $Sig$-irreducible}, \textit{regularly $Sig$-tail-irreducible}, $\bm{f} \longrightarrow_{\bm{G},\ast}\bm{h}$ are clear without any further explication. In some cases it is convenient to speak of ordinarily reducing a module element $\bm{f}\in P^m$ (i.e., without above constraint \textit{(ii)} regarding the signatures), when, in fact, we mean reducing the corresponding polynomial $\varphi(\bm{f})\in P$.

We say $\bm{f'}\in P^m$ is a (regular) $Sig$-normal form of $\bm{f}$ if
\(
\bm{f} \longrightarrow_{\bm{G},reg,\ast} \bm{f'}
\)
and $\bm{f'}$ is (regularly) $Sig$-irreducible. We denote by $\bm{f}\bmod\bm{G}$ the set of all $Sig$-normal forms of $\bm{f}$ with respect to $\bm{G}$, and by $\bm{f}\regmod\bm{G}$ the set of all regular $Sig$-normal forms with respect to $\bm{G}$.
For a pair of module elements $\bm{f}, \bm{g} \in P^m $ we define the $S$-pair of $\bm{f}$ and $\bm{g}$ as
\[
\text{Spair}(\bm{f,g})
\coloneqq \left( \frac{\ell}{LM(\varphi(\bm{f}))}\bm{f} , \frac{\ell}{LM(\varphi(\bm{g}))}\bm{g}\right)
\coloneqq \left( u\bm{f}, v\bm{g}\right)
\]
where $l \coloneqq \text{lcm}{(LT(\varphi(\bm{f})),LT(\varphi(\bm{g})))}$. We call $\text{Spair}(\bm{f,g})$ \textit{regular} if $\sig[u]{f} \neq \sig[v]{g}$ and \textit{singular} otherwise.

Let $F \coloneqq  \{f_1,...,f_m\}\subseteq P\setminus\{0\}$ be a set of polynomials, $I\coloneqq \langle F \rangle$ the ideal generated by $F$ and $\bm{s} \in T_m$ a module term. A set of module elements $\bm{G} \subseteq P^m\setminus\{\bm{0}\}$ is defined to be a \textit{$Sig$-Gröbner basis of $I$ up to signature $\bm{s}$} if
\[
\forall \bm{f} \in P^m : \sig{f} < \bm{s}\Longrightarrow 
\bm{f} \longrightarrow_{\bm{G},\ast} 0.
\]
The set $\bm{G}$ is called a \textit{$Sig$-Gröbner basis of $I$} if $\bm{G}$ is a $Sig$-Gröbner basis up to every $\bm{s}\in P^m$ (i.e., for all possible signatures $\bm{s}$). The dependence on the set $F$ (and thus the ideal $I$) is implicitly contained in the condition $\bm{f} \longrightarrow_{\bm{G},\ast} 0$, since for $\varphi=\varphi_F$ this implies \(\varphi(\bm{f})\longrightarrow_{\varphi(\bm{G}),\ast} 0.
\)\footnote{The notion of $Sig$-Gröbner bases is motivated by the fact that if $\bm{G}$ is a $Sig$-Gröbner basis, then $\varphi(\bm{G})$ is a Gröbner basis.}

A total order $\preceq$ on $\bm{G}$ with $\sig{f}\mid\sig{g}\Longrightarrow \bm{f}\preceq\bm{g}$, for all $\bm{f},\bm{g}\in\bm{G}$, is called a \textit{rewrite order}. We assume that all elements in $\bm{G}$ have distinct signatures, hence, the notion of a rewrite order is well-defined. For $\bm{s}\in T_m$, $\bm{f}\in P^m$ and $u\in T$, the element $u\bm{f}\in P^m$ is called \textit{the canonical rewriter of signature $\bm{s}$ with respect to $\bm{G}$} if $\bm{G}=\emptyset$ or if $\sig[u]{f} = \bm{s}$ and
\(
	\bm{f}=\max_{\preceq}\{\bm{g}\in\bm{G} : \sig{g}\mid\bm{s}\}.
\) 
Instead of this bulky denomination, we often just say ``the canonical rewriter of $\bm{s}$'', because the set $\bm{G}$ will be clear from the context.

\subsection{\texttt{M4GB} Algorithm}
\label{Sec:m4gb}

In 2017 Rusydi Makarim and Marc Stevens published a new algorithm for computing Gröbner bases called \texttt{M4GB}. The main innovation of \texttt{M4GB} is a fast polynomial reduction routine that only uses tail-reduced reductors in each reduction step. In addition, \texttt{M4GB} maintains a set of already used (tail-reduced) reductors and thus allows to reuse reductors. We describe a variant of $\texttt{M4GB}$ which is sketched in the performance section of \cite[Sec. 4.1]{MS17}. This variant outputs the same result as the original $\texttt{M4GB}$ algorithm, albeit it is considered more performant due to time savings in the update process of the set of reductors. The authors of \texttt{M4GB} call this variant a $\textit{lazy}$ variant, whereas we simply refer to this variant as \texttt{M4GB}. Here, we only describe the core ideas and those parts of \texttt{M4GB} that are relevant for our new Gröbner basis algorithm \algo{} in \Cref{Sec:new-algorithm}. In particular, we focus on the reduction of polynomials in \texttt{M4GB}. For a more detailed description of \texttt{M4GB} we refer the reader to the original article \cite{MS17}.

The \texttt{M4GB} algorithm essentially follows the basic outline of the textbook \texttt{Buchberger} algorithm \cite{B76}, which is ``Select, Reduce, Update'': selecting an S-pair, reducing it, and adding the reduced S-pair to the current basis in case it is nonzero. Whenever a nonzero reduced S-pair is added to the current basis, the set of S-pairs is updated. In \texttt{M4GB}, updating the set of S-pairs is achieved via the Gebauer-Möller criteria \cite{GM88}. This process is repeated until all S-pairs have been processed. In addition to the basic ``Select, Reduce, Update'' triad, \texttt{M4GB} is characterised by the following two distinct properties: $(a)$ it performs reductions only with tail-reduced reductors and, $(b)$ it maintains a list of already used (tail-reduced) reductors for future use. The benefit of these two properties are faster reductions because $(b)$ allows to reuse an already constructed (tail-reduced) reductor instead of re-constructing it again, while $(a)$ ensures that during a reduction no new reducible terms are introduced into the resulting polynomial.

More formally, let $G$ denote the current basis and $T_G(f)$ the set of reducible terms of $f\in P$ with respect to $G$. Assume \texttt{M4GB} reduces a term $t$ in a polynomial $p$ by an appropriate reductor $m$ and $m$ is \textit{not} tail-irreducible with respect to $G$. Then, for further reducing the result of the reduction $p-m$, all terms in
\[
	T_G(p-m) = \left( T_G(p)\cup T_G(m) \right)\setminus\{ t \}
\]
would have to be reduced modulo $G$. However, if $m$ is tail-irreducible we have by definition $T_G(m) \subseteq \{\text{LT}(m)\} = \{t\}$, hence
\[
	T_G(p-m) = T_G(p)\setminus\{ t \},
\]
and only terms in $T(p)\setminus \{t\}$ need to be reduced modulo $G$. This is the main conceptual advantage of $\texttt{M4GB}$ and its fast reduction routine.

Throughout all computations, \texttt{M4GB} maintains a set of reductors $M \supseteq G$, i.e., a set of monomial multiples of the current basis elements. All elements in $M$ have unique leading terms, which is why the current basis $G$ can be referenced only by its leading terms $L$. Nevertheless, we refer to $L$ as the intermediate (or current) basis. The original formulation of $\texttt{M4GB}$ in \cite{MS17} proactively updates the whole set $M$ in advance whenever a new basis element is generated. In contrast, the variant of \texttt{M4GB} that we describe (and that the authors of \cite{MS17} implement) updates the elements in $M$ only on-demand.\footnote{When we speak of updating the set of reductors $M$, this is conceptionally different from updating the set of S-pairs. The former one is specific to \texttt{M4GB}, while the latter one is an essential feature of all Gröbner basis algorithms.} This means, only when an element $m\in M$ is reused, the algorithm checks if it needs to be tail-reduced with respect to the elements referenced by $L$. This leads to a \textit{lazy} implementation of the update process of $M$. Although not explicitly stated in \cite{MS17}, for this lazy variant of \texttt{M4GB} the authors implicitly use the concept of \textit{generations}: the generation of a reductor $m\in M$ is the cardinality $|L|$ of the intermediate basis $L$ when $m$ was added to $M$. Keeping track of the generation has the following purpose: whenever a reductor $m\in M$ is reused during the execution of $\texttt{M4GB}$ and the generation of $m$ is equal to the current generation, then we know $m$ is tail-irreducible with respect to the current basis $L$ and it can be used reused without any further considerations. If the generation of $m$ is strictly smaller than the current generation, $m$ needs to updated.

\subsection{Signature-Based Algorithms}
\label{Sec:signature-algorithms}

In the textbook version of the \texttt{Buchberger} algorithm, many of the S-pairs will be reduced to zero, which means they do not contribute any new information to the eventual Gröbner basis. Hence, a reduction to zero is redundant work, and it would be nice to have an oracle detecting whether or not an S-pair will be reduced to zero \textit{without} having to carry out the actual reduction. There are criteria known to improve the textbook Buchberger algorithm in this regard (i.e., Buchberger's Product and Chain Criterion, realized in the Gebauer-Möller instantiation \cite{GM88} of the \texttt{Buchberger} algorithm), but still many redundant reductions to zero might occur. In the following, a change of perspective helps to establish even stronger criteria for detecting redundant reductions to zero. Let $f$ be a polynomial in the ideal generated by the polynomials $f_1,\ldots,f_m\in P$, i.e., $f \in \langle f_1,...,f_m\rangle$. Then $f$ can be written as
\(
	f = \sum_{i = 1}^m p_i f_i,
\)
for some polynomials $p_1,\ldots,p_m\in P$ (which are not necessarily unique). This notation of $f$ motivates a new perspective: $f$ cannot only be considered as polynomial but also as module element $(p_1,\ldots,p_m)\in P^m$. Adopting the module's perspective, it is possible to introduce a new concept called \textit{signatures} for detecting unnecessary S-pair reductions. The main idea behind signatures is, roughly speaking, to keep track of how the polynomials generated during a Gröbner basis computation depend on the original input polynomials. More concretely, this means a signature-based algorithm not only processes information coming from a polynomial $f$ itself but also from the vector $(p_1,\ldots,p_m)$ constituting the relation $f=\sum_i p_i f_i$, where the $f_i$ would be the original input polynomials. On the one hand, this idea aims at exploiting zero-relations between the input polynomials (i.e., syzygies from the module perspective) to detect redundant reductions; on the other hand, it uses the (more subtle) fact that different polynomial combinations of the input polynomials (i.e., different module elements from the module perspective) can have the same reduction remainder. Thus only one of these reductions need to be performed. 
The former observation is the basis for the so-called \textit{syzygy criterion}, while the latter observation leads to the \textit{rewrite-criterion} (see \Cref{Line:new-syzygy} and \ref{Line:rewriter}, respectively, in \Cref{Alg:main}).

With above motivation of signatures at hand, we state the signature equivalent of Buchberger's S-pair criterion. The fundamental theorem underlying all signature-based algorithms is the following result.
\begin{theorem}[\cite{ER13}, Theorem 3]
	\label{Thm:sig-S-pair-criterion}
	Let $\bm{s}\in T_m$ be a module term and $\bm{G} \subseteq P^m$ be a finite set of module elements. If for all $\bm{p}\in P^m$ with $\bm{p}$ a regular S-pair of elements in $\bm{G}$ or $\bm{p}$ a canonical basis vector $\bm{e_i}$ 
	(and $\sig{p} \bm{<} \bm{s}$, resp.) it holds that $\bm{p}\regmod\bm{G}$ contains a syzygy or a singularly $Sig$-top-reducible element, then $\bm{G}$ is a signature Gröbner basis (up to $\bm{s}$, resp.).
\end{theorem}
The following two observations explicate how signatures help to detect unnecessary reductions to zero in advance: assume we have a Gröbner basis $\bm{G}\subseteq P^m$ up to signature $\bm{s}\in T_m$.
First, one can show that for any two regularly $Sig$-irreducible module elements $\bm{f},\bm{g}\in P^m$ it holds
\[
	\sig{f} = \sig{g} = \bm{s}\Longrightarrow \varphi(\bm{f}) = c\cdot \varphi(\bm{g})
\]
for some $c\in\mathbb{F}\setminus \{0\}$. Second, if there exists a syzygy $\bm{h}\in P^m$ with $\sig{h}\mid \bm{s}$, then
\[
\forall \bm{f}\in P^m\text{ with }\sig{f}=\bm{s}: \bm{f}\longrightarrow_{\bm{G},\ast} 0.
\]
The salient points are:
\textit{(a)} we only need to $Sig$-reduce one element with a given signature (we will choose the one which is `easier' to handle). Hence, in a signature-based algorithm, instead of an S-pair with a given signature, we are free to choose any module element with the same signature and reduce this element to check whether the current signature provides new information for our eventual Gröbner basis. This approach is called \textit{rewriting} and Gröbner basis algorithms based on this approach are called \textit{rewrite} algorithms \cite{ER13}; \textit{(b)} if we know that the signature of the element to be reduced is a multiple of the signature of a syzygy, we can skip the computation of the reduction at all. This is why a signature-based algorithm always keeps track of syzygy signatures and stores them separately.

This is all we intend to say about the ideas behind signature-based and rewrite Gröbner basis algorithms and, in particular, we do not state a pseudo code for them. The basic ideas we adopt from the signature and rewriting approach for our \algo{} algorithm are evident from Algorithm \ref{Alg:main}.
For a more in-depth motivation and treatment of signature-based and rewrite Gröbner basis algorithms we refer to the comprehensive survey article \cite{EF17}.

\section{\algo{} Algorithm}
\label{Sec:new-algorithm}

In this section, we present our new Gröbner basis algorithm \algo{} that amalgamates the core ideas of (signature-based) rewrite algorithms with the main ideas of \texttt{M4GB}. For this amalgamation to be viable, we introduce a new concept called \textit{signature flags}. On a high level, signature flags play a similar role as generations in \texttt{M4GB} and allow to \textit{efficiently} fuse the ideas behind signature-based algorithms and \texttt{M4GB}, respectively. As such, \algo{} is an algorithm which aims to combine the strengths of both worlds:
\textit{(a)} fast reduction of polynomials due to the \texttt{M4GB}-like reduction routine;
\textit{(b)} strong criteria for discarding redundant S-pairs adopted from signature-based algorithms.

\subsection{New Definitions}
\label{Sec:new-definitions}

Since \algo{} works with $Sig$-tail-irreducible reductors up to some signature $\bm{s}$, we explicate this concept in a formal definition. In the following let $\bm{G}\subseteq P^m\setminus\{\bm{0}\}$ be a non-empty and finite set of non-zero module elements.

We call a term $t\in T$ \textit{$Sig$-reducible with respect to $\bm{G}$ and up to $\bm{s}$}, if there exist $u\in T$, $\bm{g}\in\bm{G}$ such that $\text{LT}(\varphi(u\bm{g}))=t$ and $\sig[u]{g}< \bm{s}$. A module element $\bm{f}\in P^m$ is called \textit{$Sig$-reducible with respect to $\bm{G}$ and up to $\bm{s}$}, if there exists a term $t\in T(\varphi(\bm{f}))$ that is $Sig$-reducible with respect to $\bm{G}$ and up to $\bm{s}$. We denote such a reduction step by $\bm{f}\longrightarrow_{\bm{G},\bm{s}}\bm{f}-c\cdot u\bm{g}$, for an appropriate scalar $c\in\mathbb{F}\setminus\{0\}$. For a given set of terms $D\subseteq T(\varphi(\bm{f}))$, we call $\bm{f}$ \textit{$Sig$-reducible with respect to $\bm{G}$, $D$ and up to $\bm{s}$} if there exists a reductor $u\bm{g}$ of $\bm{f}$ such that $\text{LT}(\varphi(u\bm{g}))=d$ for some $d\in D$ and $\sig[u]{g}< \bm{s}$. Such a reduction step is denoted by $\bm{f}\longrightarrow_{\bm{G},\bm{s},D}\bm{f}-c\cdot u\bm{g}$.

\begin{remark}
In particular, for $\bm{s}=\sig{f}$ and $D=T(\varphi(\bm{f}))$, the reduction $\bm{f}\longrightarrow_{\bm{G},\bm{s},D}\bm{f}-c\cdot u\bm{g}$ describes a regular $Sig$-reduction $\bm{f}\longrightarrow_{\bm{G},reg}\bm{f}-c\cdot u\bm{g}$. This means, our new view $\longrightarrow_{\bm{G},\bm{s},D}$ on signature-based reductions contains regular $Sig$-reductions as special case. Moreover, if we choose $D=T(\tail{\varphi(\bm{f})})$, we allow all regular $Sig$-reductions except for a top-reduction. These two special cases are the instantations of $D$ we are most interested in, although the statements below, e.g., Lemma \ref{Lem:correct-reduce}, can be applied for arbitrary $D \subseteq T(\varphi(\bm{f}))$.
\end{remark}

As highlighted in \Cref{Sec:preliminary-definitions}, we often do not explicitly mention the set $\bm{G}$. In the same manner, we define $\bm{f}\longrightarrow_{\bm{G},\bm{s},\ast}\bm{h}$, $\bm{f}\longrightarrow_{\bm{G},\bm{s},D,\ast}\bm{h}$, \textit{$Sig$-irreducible up to $\bm{s}$}, \textit{$Sig$-irreducible with respect to $D$ and up to $\bm{s}$}, \textit{$Sig$-tail-irreducible up to $\bm{s}$}. A \textit{normal form of $\bm{f}$ with respect to $\bm{G}$ and up to $\bm{s}$} is an element $\bm{f'}\in P^m$ that is $Sig$-irreducible with respect to $\bm{G}$ and up to $\bm{s}$ and for which it holds $\bm{f}\longrightarrow_{\bm{G},\bm{s},\ast} \bm{f'}$. We denote the set of all normal forms of $\bm{f}$ with respect to $\bm{G}$ and up to $\bm{s}$ by $\bm{f}\bmod\!_{\!\bm{s}}\;\bm{G}$. For a set of terms $D\subseteq T(\varphi(\bm{f}))$, a \textit{normal form of $\bm{f}$ with respect to $\bm{G}$, $D$ and up to $\bm{s}$} is an element $\bm{f'}\in P^m$ that is $Sig$-irreducible with respect to $\bm{G}$, $D$ and up to $\bm{s}$ such that $\bm{f}\longrightarrow_{\bm{G},\bm{s},D,\ast} \bm{f'}$.  We denote the set of all normal forms of $\bm{f}$ with respect to $\bm{G}$ and up to $\bm{s}$ by $\bm{f}\bmod\!_{\!\bm{s},D}\;\bm{G}$.

As in \texttt{M4GB}, the \textit{generation} $\gen{\bm{m}}$ of a reductor $\bm{m}\in\bm{M}$ is defined as the cardinality of the set $\bm{G}$ at the time $\bm{m}$ is constructed. \footnote{Here, the term ``constructed'' also encompasses the case when $\bm{m}$ is updated, or in other words, ``re-constructed''.} We denote the instance of $\bm{G}$ at this time with $\bm{G}_{\gen{\bm{m}}}$.
For a module element $\bm{f}\in P^m$ the \textit{signature flag with respect to $\bm{G}$} is defined as
\[
	\flag{\bm{f}}\coloneqq  \min\{
		\sig[v]{g} : \bm{g} \in \bm{G}, v \in T, v\bm{g}\textit{ a tail-reductor of }\bm{f}
	\},
\]
or $\flag{\bm{f}}\coloneqq \bm{\infty}$ if $\bm{f}$ is tail-irreducible.
The symbol $\bm{\infty}$ can be understood as a formal symbol added to $T_m$ with the simple property that
\[\forall \bm{s} \in T_m: \bm{s} < \bm{\infty}.\]

\subsection{Description of \algo{}}
\label{Sec:description-algorithm}

The overall structure of \algo{} is depicted in Algorithm \ref{Alg:main} and resembles the basic structure of a rewrite Gröbner basis algorithm (as outlined in \cite{ER13}) with signature-based criteria to discard redundant S-pairs (see Line \ref{Line:rewriter}, \ref{Line:same-signature}, \ref{Line:syzygy}) and the fundamental ``Select, Reduce, Update'' triad from the \texttt{Buchberger} algorithm \cite{B76}. In particular, \algo{} processes S-pairs in strictly increasing signature and keeps track of syzygy signatures in a separate set $\bm{H}$ (Line \ref{Line:new-syzygy}). If a new basis element is found (Line \ref{Line:new-element}), the \updatealgo{} routine (\Cref{Alg:update}) for the current basis $\bm{G}$ and the current set of S-pairs $\bm{P}$ is triggered. The steps in \updatealgo{} are governed by the same principles as in any other signature-based algorithm, with the difference, that \updatealgo{} detects whether a basis element $\bm{e_i}$ has been processed and thus extends the set of syzygy signatures $\bm{H}$ accordingly (Algorithm \ref{Alg:update}, Line \ref{Line:extend-syzygy-signatures}). The main innovations of \algo{} are incorporated into the reduction routine \reducealgo{} described in Algorithm \ref{Alg:reduce}. In the following, we discuss the novel features as well as the intricacies of \reducealgo{} more comprehensively.

\begin{algorithm}[H]
    \caption{\algo{}}
    \label{Alg:main}
	\SetKwProg{Ifi}{if}{ then}{}
	\SetKwProg{Whilei}{while}{}{}
	\SetKwBlock{Begin}{if}{}
	\KwIn{Non-zero input polynomials $F = \{f_1,...,f_m\}$, a rewrite order, a term order on $T$ and a compatible order extension on $T_m$}
	\KwOut{A Gröbner basis $G$ of the ideal generated by $F$}
	$\bm{G} \coloneqq \emptyset$;	$\bm{M} \coloneqq \emptyset$;  
	$\bm{H} \coloneqq \emptyset$ \\
	$\bm{P} \coloneqq \{\bm{e_i}: i \in \{1,...,m\}\}$\\
 	\Whilei{$\bm{P}\neq \emptyset$} {
 		Select $\bm{f}\in\bm{P}$ with minimal signature $\bm{s}=\sig{f}$ and $u\bm{g}$ the canonical rewriter of $\bm{s}$ w.r.t $\bm{G}$.\label{Line:rewriter}\\
 		$\bm{P} \coloneqq \bm{P} \!\setminus\!\{\bm{p}\in\bm{P} : \sig{p} = \bm{s} \}$\label{Line:same-signature}\\
		\If{$\bm{s}$ is not divisible by some $\bm{h}\in\bm{H}$\label{Line:syzygy}}{
			($\bm{M}$, $\bm{f'}$) $\coloneqq  \reducealgo{}(\bm{f},\varphi(\bm{f}),\bm{s},\bm{M},\bm{G})$\\
 			\If{$\bm{f'} = 0$\label{Line:new-syzygy}}{
				$\bm{H} \coloneqq \bm{H} \cup \{\bm{s}\}$\\
 			}
  		    \Else{\label{Line:new-element}
                ($\bm{G}$, $\bm{P}$, $\bm{H}$) $\coloneqq$ 
                \updatealgo($\bm{f'}$,  $\bm{G}$, $\bm{P}$, $\bm{H}$)\\
			}
    	} 
	}
	\Return $\varphi(\bm{G})$\\
\end{algorithm}

\begin{algorithm}[H]
	\caption{\texttt{Update}}
	\label{Alg:update}
	\KwIn{Current basis $\bm{G}$, set of S-pairs $\bm{P}$ and set of syzygy signatures $\bm{H}$, new basis element $\bm{f}$}
	\KwOut{Updated $\bm{G}$, $\bm{P}$ and $\bm{H}$}
    \vspace{1em}
	$\bm{P}\coloneqq  \bm{P}\cup\{\text{Spair}(\bm{f},\bm{g}) : \bm{g}\in\bm{G},\,\text{Spair}(\bm{f},\bm{g})\text{ regular} \}$\\
    $\bm{G}\coloneqq \bm{G}\cup \{ \bm{f} \}$\\
	\If(\tcp*[f]{\normalfont $\bm{f}$ comes from a basis element $\bm{e_i}$})
	{$\sig{f}$ = $\bm{e_i}$\label{Line:extend-syzygy-signatures}}{
	    $\bm{H} \coloneqq  \bm{H} \cup \{\sig[\varphi(\bm{g})\bm{e_i}-\varphi(\bm{e_i})\bm{g}]{}: \bm{g} \in \bm{G}\}$\\
	}
    \Return ($\bm{G}$, $\bm{P}$, $\bm{H}$)
\end{algorithm}

As in \texttt{M4GB}, the \reducealgo{} routine keeps track of previously used reductors and stores them in a set $\bm{M}$. The key feature of \texttt{M4GB}, namely, working with tail-reduced reductors, is implemented in \reducealgo{} as well. The difference to \texttt{M4GB} and a crucial point is that whenever a reductor $\bm{m}\in P^m$ is added to $\bm{M}$, it need not be fully tail-irreducible with respect to $\bm{G}$ but only $Sig$-tail-irreducible up to the current signature $\bm{s}$. This property is an important part of our efficient amalgamation of signature-based algorithms with \texttt{M4GB}: by \Cref{Thm:sig-S-pair-criterion}, signature-based algorithms work with \textit{regular $Sig$-reductions} and hence, only those terms in $T(\tail{\bm{m})})$ need be reduced that have a reducer with signature smaller than $\bm{s}$. We formalized this particular property in the definitions in \Cref{Sec:new-definitions}.

Again, as in \texttt{M4GB}, elements in $\bm{M}$ are updated in a lazy manner, meaning only on-demand when they are reused and not proactively whenever a new basis element is added to $\bm{G}$.
\begin{remark}
In the context of \algo{}, updating an element of $\bm{M}$ alludes to the process of restoring its $Sig$-tail-irreducibility with respect to the current basis $\bm{G}$ and up to the current signature $\bm{s}$.
\end{remark}
Below, we consider the two scenarios when some element $\bm{m}\in\bm{M}$ stops being $Sig$-tail-irreducible and thus needed to be updated in case it was reused:
\begin{enumerate}[\itshape(1)]
\item\label{Item:new-basis-element}
 A new element is added to $\bm{G}$ which regularly $Sig$-tail-reduces $\bm{m}$.
\item\label{Item:existing-tail-reductor}
An existing tail-reductor of $\bm{m}$ becomes a valid $Sig$-tail-reductor in light of the current signature $\bm{s}$.
\end{enumerate}
The aspect in \textit{(\ref{Item:existing-tail-reductor})} needs some clarification. Assume, at the time $\bm{m}$ was added to $\bm{M}$ it was regularly $Sig$-tail-irreducible up to some signature $\bm{r}$ but not ordinarily tail-irreducible (i.e., $\varphi(\bm{m})$ is not tail-irreducible with respect to $\varphi(\bm{G})$).
This means, at the time $\bm{m}$ was added to $\bm{M}$ there was some basis element multiple $u\bm{g}$ which tail-reduced $\bm{m}$ but the reduction was not a valid regular $Sig$-tail-reduction up to $\bm{r}$, because $\bm{r}\leq\sig[u]{g}$. If the current signature $\bm{s}$ fulfills $\bm{s} > \sig[u]{g}$, the reductor $u\bm{g}$ becomes a valid reductor for a regular $Sig$-tail-reduction.

\begin{algorithm}[t]
	\caption{\reducealgo{}}
	\label{Alg:reduce}
	\KwIn{$\bm{f}\in P^{m}$, polynomial $p\in P$ with $T(p)\subseteq T(\varphi(\bm{f}))$, signature $\bm{s}$, current basis $\bm{G}$, current set of reductors $\bm{M} \subseteq P^m$}
	\KwOut{Possibly extended set $\bm{M}$, $Sig$-normal form $\bm{f'}\in \bm{f}\bmod\!_{\!\bm{s},T(p)}\;\bm{G}$ with respect to $T(p)$ and up to $\bm{s}$}
	$\bm{f'} \coloneqq  \bm{f}$\\
	\For{$t\in T(p)$\label{Line:main-loop}}
	{	
		\If{$\exists\bm{m}\in\bm{M} : \text{LT}(\varphi(\bm{m}))=t$\label{Line:reductor-present}}{
			Select such $\bm{m}$\\
			$\bm{m'}\coloneqq \bm{m}$\\
			\If{$\gen{\bm{m}}<|\bm{G}|$\label{Line:smaller-generation}}{
				$\bm{M} \coloneqq  \bm{M} \!\setminus\! \{ \bm{m} \}$\\
				$(\bm{M}, \bm{m'}) \coloneqq  
				\reducealgo(\bm{m}, \text{Tail}(\varphi(\bm{m})), \bm{s},\bm{G}\backslash\bm{G}_{\gen{\bm{m}}},\bm{M})$\label{Line:update-because-new-basis-element}\\
				\If{$Flag(\bm{m}) < \bm{s}$\label{Line:smaller-flag-1}}{
					$(\bm{M}, \bm{m'}) \coloneqq  
					\reducealgo(\bm{m'},\text{Tail}(\varphi(\bm{m'})),\bm{s},\bm{G}_{\gen{\bm{m}}},\bm{M})$\label{Line:old-gen-old-flag}\\
				}
				$(\bm{M},\bm{m'}) \coloneqq  \texttt{UpdateM}(\bm{M},\bm{m'},\bm{G})$\\
			}
			\ElseIf{$Flag(\bm{m})< \bm{s}$\label{Line:smaller-flag-2}}{
				$\bm{M} \coloneqq  \bm{M} \!\setminus\! \{ \bm{m} \}$\\
				$(\bm{M},\bm{m'}) \coloneqq  \reducealgo(\bm{m'},\text{Tail}(\varphi(\bm{m'})),\bm{s},\bm{G},\bm{M})$\label{Line:same-gen-smaller-flag-reduce}\\
				$(\bm{M},\bm{m'}) \coloneqq  \texttt{UpdateM}(\bm{M},\bm{m'},\bm{G})$\\
			}
			$\bm{f'} \coloneqq  \bm{f'} - C_t(\varphi(\bm{f'}))\cdot\bm{m'}$\\
		}
		\ElseIf{$\exists \bm{g}\in\bm{G}: \text{LT}(\varphi(\bm{g}))\mid t$, $\sig[u]{g}<\bm{s}$, $u\coloneqq t\slash\text{LT}(\varphi(\bm{g}))$\label{Line:reductor-not-present}}
		{
			
			Select such $\bm{g}$\\ 
			$(\bm{M}, \bm{m'}) \coloneqq  \reducealgo(u\bm{g},\text{Tail}(\varphi(u\bm{g})), \bm{s}, \bm{G},\bm{M})$\\
			$(\bm{M},\bm{m'}) \coloneqq  \texttt{UpdateM}(\bm{M},\bm{m'},\bm{G})$\\
			$\bm{f'} \coloneqq  \bm{f'} - C_t(\varphi(\bm{f'}))\cdot\bm{m'}$\\
			
		}
	}
	\Return ($\bm{M}$, $\bm{f'}$);
\end{algorithm}

To resolve \textit{(\ref{Item:new-basis-element})}, we use the concept of `generations' (adopted from \texttt{M4GB}). All reductors added to $\bm{M}$ are equipped with a generation (the cardinality of $\bm{G}$ at the time $\bm{m}$ is created). Everytime a new basis element is added to $\bm{G}$, the generation increases and thus, any reductor in $\bm{M}$ being reused and having a strictly smaller generation than the current one needs to be updated (Algorithm \ref{Alg:reduce}, Line \ref{Line:smaller-generation}).
To resolve \textit{(\ref{Item:existing-tail-reductor})}, we use the new concept of `signature flags'. All reductors added to $\bm{M}$ are equipped with a signature flag. The idea of a signature flag is to define it as the minimal signature for which \textit{(\ref{Item:existing-tail-reductor})} occurs. Consequently, if the signature flag of a reductor being reused is smaller than the current signature, the reductor needs to be updated (Algorithm \ref{Alg:reduce}, Line \ref{Line:smaller-flag-1} and \ref{Line:smaller-flag-2}).

\begin{algorithm}[t]
	\caption{\texttt{UpdateM}}

	\KwIn{Current set of reductors $\bm{M}$, reductor $\bm{m'}$ to be normalized and equipped with generation and signature flag, current basis $\bm{G}$}
	\KwOut{Updated set $\bm{M}$ and updated $\bm{m'}$}

	$\bm{m'} \coloneqq  LC(\varphi(\bm{m'}))^{-1}\cdot\bm{m'}$\\
	$Flag(\bm{m'}) \coloneqq  \min\{Flag(t): t \in T(\tail{\varphi(\bm{m'})})\}$\\
	$\gen{\bm{m'}}\coloneqq |\bm{G}|$\\
	$\bm{M} \coloneqq  \bm{M} \cup \{\bm{m'}\}$\\

	\Return $(\bm{M},\bm{m'})$
\end{algorithm}

\subsection{Termination and Correctness}
Before we prove termination and correctness, we want to shed more light on the particular update process of reductors in \reducealgo{}. For this, we come back to the two situations in \Cref{Sec:description-algorithm} when a reductor $\bm{m}\in\bm{M}$ stops being $Sig$-tail-irreducible with respect to the current basis $\bm{G}$ and up to the current signature $\bm{s}$. Here, we state them more formally and by means of our new definitions from \Cref{Sec:new-definitions}. Case \textit{(\ref{Item:new-basis-element})} in \Cref{Sec:description-algorithm} corresponds to
\[
	\exists t \in T(\tail{\bm{m}}), \bm{g} \in \bm{G}\setminus\bm{G}_{ \gen{\bm{m}} }, v \in T: \sig[v]{g} < \bm{s}\wedge\text{LT}(\varphi(v\bm{g})) = t,
\]
whereas case \textit{(\ref{Item:existing-tail-reductor})} is characterised by
\[
	\exists t \in T(\tail{\bm{m}}), \bm{g} \in \bm{G}_{\gen{\bm{m}}}, v \in T: \sig[v]{g}<\bm{s}\wedge \text{LT}(\varphi(v\bm{g})) = t.
\]
This is the reason why \reducealgo{} only needs to $Sig$-reduce with respect to $\bm{G}\setminus\bm{G}_{\gen{\bm{m}}}$ in Line \ref{Line:update-because-new-basis-element} whenever an update due to an older generation is necessary and the same reasoning applies to Line \ref{Line:old-gen-old-flag} and $\bm{G}_{\gen{\bm{m}}}$.

The outline of \algo{} follows the same outline as a rewrite basis algorithm, with only the reduction routine \reducealgo{} being different. Since \algo{} always calls \reducealgo{} with the arguments $(\bm{f},\varphi(\bm{f}), \bm{s},\bm{M},\bm{G})$ and it holds $\bm{s}=\sig{f}$,
we only need to prove correctness and termination of \reducealgo{} to argue correctness and termination for \algo{}. We begin with an important lemma. In essence, \Cref{Lem:correct-reduce} explains why \reducealgo{} correctly computes a $Sig$-normal form with respect to a given set of terms $D$ and up to signature $\bm{s}$. We emphasize that the usage of $Sig$-tail-irreducible reductors is crucial here, without it, the statement would be wrong.

\begin{lemma}
\label{Lem:correct-reduce}
Let $\bm{f}\in P^m$, $\bm{s}\in T_m\cup\{\bm{\infty}\}$ and $\bm{G}\subseteq P^m$. Let $T_{\bm{G},\bm{s},D}(\bm{f})$ denote the set of all terms in $D\subseteq T(\varphi(\bm{f}))$ that are $Sig$-reducible with respect to $\bm{G}$ and up to signature $\bm{s}$. For each $t\in T_{\bm{G},\bm{s},D}(\bm{f})$, let $\bm{m}_t$ denote a reductor of $t$ with $\sig{m_t}<\bm{s}$ which is $Sig$-tail-irreducible with respect to $\bm{G}$ and up to $\bm{s}$. Then
\[
\bm{f'}\coloneqq \bm{f}-\sum_{t\in T_{\bm{G},\bm{s},D}(\bm{f})} c_t(f)\cdot\bm{m_t} \in \bm{f}\bmod\!_{\!\bm{s},D} 
\;\bm{G}.
\]
\end{lemma}

\begin{proof} Because all $\bm{m}_t$ are $Sig$-tail-irreducible up to $\bm{s}$, we have
\begin{align*}
	T_{\bm{G},\bm{s},D}(\bm{f'})
	&\subseteq \Big(T_{\bm{G},\bm{s},D}(\bm{f})\cup\bigcup_{t\in T_{\bm{G},\bm{s},D}(\bm{f})}T_{\bm{G},\bm{s},D}(\bm{m}_t)\Big)\setminus T_{\bm{G},\bm{s},D}(\bm{f})
	&= \emptyset,
\end{align*}
so it follows that $\bm{f'}$ is $Sig$-irreducible with respect to $D$ and up to $\bm{s}$. We are left to show that $\bm{f}\longrightarrow_{\bm{G},\bm{s},D,\ast}\bm{f'}$.
To do so, we proceed inductively: assume by hypothesis that for a fixed $n \in \mathbb{N}$, we have that
\[\bm{f} \longrightarrow_{\bm{G},\bm{s},D,\ast} \bm{f_n} := \bm{f}-\sum_{t\in S} c_t(f)\cdot\bm{m_t} \]
holds for arbitrary $S \subseteq T_{\bm{G},\bm{s},D}(\bm{f})$ with $\lvert S \rvert = n$.
 If $S = T_{\bm{G},\bm{s},D}(\bm{f})$, then $\bm{f'}=\bm{f_n}$ and the claim holds trivially.
 Otherwise, let $t_0 \in T_{\bm{G},\bm{s},D}(\bm{f})\setminus S$.
 We need to show that $\bm{f_n}\longrightarrow_{\bm{G},\bm{s},D,\ast} \bm{f_n}- c_{t_0}\bm{m_{t_0}}$.
 As $\sig{\bm{f}} = \sig{\bm{f_n}}$, it suffices to show that $t_0 \in T_{\bm{G},\bm{s},D}(\bm{f_n})$.
  As $\bm{m_t}$ is Sig-tail-irreducible by assumption for
 every $t \in S$, we have $t_0 \notin \bigcup_{t \in S}T(\bm{m_t})$ and in particular,
 $t_0 \in T_{\bm{G},\bm{s},D}(\bm{f_n})$ follows. This concludes the proof.
\end{proof}

\begin{theorem} \reducealgo{} terminates and correctly computes a $Sig$-normal form $\bm{f'}\in \bm{f}\bmod\!_{\!\bm{s},T(p)}\;\bm{G}$.
\end{theorem}

\begin{proof}
For termination, we note that whenever \reducealgo{} is processing a term $t$ in recursion level $n\in\mathbb{N}_0$ and calls itself, all terms being processed in the following recursion level $n+1$ regarding $t$ are strictly smaller than $t$. This is because whenever \reducealgo{} calls itself in level $n$ while processing a term $t$, it calls itself on $\tail{v}$ of some polynomial $v$ with $\text{LT}(v)=t$ and thus for any subsequent term $u$ in level $n+1$ regarding $t$ it holds $u<t$. Hence, the recursion depth of \reducealgo{} must be finite. Since at a given recursion level only finitely many terms are being processed, we conclude that \reducealgo{} eventually terminates.

To argue correctness, in view of \Cref{Lem:correct-reduce}, it suffices to prove that for every $t\in T(p)$, a potential reductor $\bm{m'_t}$ is $Sig$-tail-irreducible up to $\bm{s}$ and fulfills $\sig{m'_t}<\bm{s}$.

It is clear that \reducealgo{} reaches the end of a recursive path if and only if it processes a reductor where it does not call itself anymore. Looking at Algorithm \ref{Alg:reduce}, this is the case if and only if \reducealgo{} is being called with $(\bm{f},T(p),\bm{s},\bm{G},\bm{M})$ such that every $t\in T(p)$ is either \textit{(i)} $Sig$-irreducible with respect to $\bm{G}$ and up to $\bm{s}$ or \textit{(ii)} there already exists a reductor $\bm{m}\in\bm{M}$ with $\gen{\bm{m}}=|\bm{G}|$ and $\text{Flag}(\bm{m})\geq\bm{s}$. By the definitions of generation and signature flag and by the construction of elements in $\bm{M}$, \textit{(ii)} is equivalent with $\sig{m}<\bm{s}$ and $\bm{m}$ being $Sig$-tail-irreducible with respect to $\bm{G}$ and up to $\bm{s}$. Using \Cref{Lem:correct-reduce}, we deduce that the reduction remainder $\bm{f'} \in \bm{f}\bmod\!_{\!\bm{s},T(p)} 
\;\bm{G}.$ If $\bm{f'}$ serves as a reductor $\bm{m'_t}$ in a recursion level above, note that $\sig{\bm{f}} < \bm{s}$ and hence,
also $\sig{\bm{f'}} < \bm{s}$ follows.
\end{proof}

\section{Implementation \& Performance}

In this section, we discuss some implementation details and the performance of our \algo{} algorithm. We base our implementation on the Mathic \verb|C++|-library developed by Roune \cite{ML}. In Mathic, we integrate our algorithm as a new module into the MathicGB Gröbner basis module. Using the Mathic framework allows us to directly compare the performance of \algo{} against the signature-based algorithm \texttt{SB} presented by Roune and Stillman \cite{RS12}. Keeping the same naming convention as in \cite{RS12}, we refer to their signature-based Gröbner basis algorithm as \texttt{SB}. As we optimize our implementation of the reduction routine outlined in \Cref{Alg:reduce}, there are minor differences to the pseudocode. These differences have no impact on the overall behaviour or correctness of the algorithm. Instead, they aim to leverage the language-specific advantages of \texttt{C++} to create a competitive proof-of-concept implementation. The source code of our implementation of \algo{} is available under
\url{https://extgit.iaik.tugraz.at/krypto/m5gb.git}.

We show that using the same library for implementing \texttt{SB} and \algo{}, we obtain a significant, scalable speed-up for dense, quadratic, overdefined polynomial systems. These systems are used for benchmark purposes in the original article about \texttt{M4GB} by Makarim and Stevens \cite{MS17} and are posed as a problem instance in the MQ Challenge \cite{mq}.

We also performed informal tests for other systems, e.g., some canonical test systems in the literature like \texttt{katsura}, \texttt{eco} or \texttt{cyclic}. Most of the results indicated that the performance of \algo{} falls behind that of \texttt{SB}. We conjecture several reasons behind these results. First, creating tail-reduced reductors is time-consuming and, depending on the structure of the polynomial system, may not yield an overall advantage compared to using ordinary reductors. Second, due to the recursive nature of the \texttt{M4GB}-style reductions, we cannot use the efficient data structures that Mathic uses to increase the performance of their implementation. Lastly, and connected to the previous point, since \algo{} uses \texttt{M4GB}-style reduction, our algorithm also inherits the disadvantages of \texttt{M4GB}. This is further evidenced by the outcomes of informal comparisons between \texttt{M4GB} and \algo{}. Although these two algorithms are implemented in a substantially different way, we found that whenever \algo{} performed poorly this also was the case for \texttt{M4GB}. However, to provide a more reliable conclusion in this regard, further and more systematic experiments are needed. This, as well, includes implementing \texttt{M4GB} and \algo{} in a more comparable manner. We leave this open for future work. 

\subsection{Implementation Details}
The original signature-based algorithm \texttt{SB} of Roune and Stillman \cite{RS12,mathicgb} does not use signature flags and generations. Thus, we extend the underlying data structures such that generations and signature flags are supported. Both generations and signature flags are implemented on term and polynomial granularity. Each polynomial stores its generation as an integer value. An unordered map $\bm{I}$, that maps term hashes to generations, stores the generations of irreducible terms. We do not need to store additional information for reducible terms, as they always cause a reducer lookup in the current basis or a lookup in the current set of tail-reduced reductors $\bm{M}$.  

Contrary to the pseudocode, we do not explicitly calculate and store signature flags when terms and module elements are stored in $\bm{I}$ or $\bm{M}$. Instead, we only store the information on whether a term has a dividing leading term in the base. Only those terms may have a signature flag that is not $\bm{\infty}$. We encode this information using a single bit in the generation integer. If a term or polynomial has a finite flag, we calculate and store the actual signature flag on the first subsequent access. This approach allows us to calculate signature flags only for elements where the flag is actually needed by the algorithm. Thus, we avoid unnecessary flag computations and also reduce the memory overhead. For systems where most signature flags are infinite, this optimization allows us to skip most of the flag logic, which leads to a further increase in the performance of our implementation.
 
\subsection{Performance Metrics}
We evaluate our implementation of \algo{} by computing the Gröbner basis for overdefined dense quadratic systems with an increasing number of variables $N$ and $M=2N$ polynomials over $\mathbb{F}_{101}$, unless stated otherwise. The variable count ranges from $N=5$ to $N=21$. For each $N$, we generate 10 distinct equation systems that are certain to have a solution. The performance metrics for each $N$ are computed as the arithmetic mean of the metric over all 10 system instances. In our evaluation, we consider the following three metrics.

\paragraph*{Time per Basis Element} The time spent per basis element is the primary indicator of the performance of our algorithm. A lower amount of time per basis element indicates a faster implementation. The resulting time per basis element is computed as the overall runtime divided by the number of elements that reside in the final Gröbner basis.

\paragraph*{Peak Memory Usage} We monitor the memory consumption of the implementations using the \texttt{time}-program on Linux. While we could track all memory allocations in the program through instrumentation-based monitoring, we chose to measure the overall memory footprint instead, as it can become a limiting factor when calculating large bases.

\paragraph*{Number of Reductions} As a third metric, we keep track of the number of actual reductions. A \textit{reduction} (or \textit{reduction step}) in this context is a single step in the process of reducing a polynomial (with respect to some set of divisors). We extend the existing \texttt{SB} implementation such that each reduction step is counted. Likewise, we keep track of reduction steps in \algo{} as well. Then, for a fixed polynomial system, we compare the respective number of reductions in \texttt{SB} and \algo{}.

\subsection{Evaluation and Discussion}

\begin{figure}[t]
	\begin{subfigure}[t]{0.49\textwidth}
	  \begin{tikzpicture}
    \pgfplotstabletranspose[input colnames to={n},colnames from={n}]{\mytabletoplot}{data/timing.csv}
    \begin{axis}[
        grid=both,
        major grid style = {lightgray},
        minor grid style = {lightgray!25},
        width=\textwidth,
        ymode = log,
        ylabel={Time in \SI{}{\milli\second}},
        xlabel={\# of Variables $N$},
        legend pos= north west,
        ]
        \addplot+[] table [x = {n}, y = {SigGb}] {\mytabletoplot};
        \addlegendentry{\texttt{SB}};
        \addplot+[] table [x = {n}, y = {Chimera}] {\mytabletoplot};
        \addlegendentry{\algo};
    \end{axis}
\end{tikzpicture}
		\caption{
		The mean time spent per basis element calculated over 10 dense quadratic systems per variable count.
	}
	\label{fig:dq_timing1}
	\end{subfigure}
	\hfill
	\begin{subfigure}[t]{0.49\textwidth}
		\begin{tikzpicture}
    \pgfplotstabletranspose[input colnames to={n},colnames from={n}]{\mytabletoplot}{data/mem.csv}
    \begin{axis}[
        grid=both,
        major grid style = {lightgray},
        minor grid style = {lightgray!25},
        width=\textwidth,
        ymode = log,
        legend pos= north west,
        ylabel={Memory in \SI{}{\kilo\byte}},
        xlabel={\# of Variables $N$},
        ]
        \addplot+[] table [x = {n}, y = {SigGb}] {\mytabletoplot};
        \addlegendentry{\texttt{SB}};
        \addplot+[] table [x = {n}, y = {Chimera}] {\mytabletoplot};
        \addlegendentry{\algo};
       
    \end{axis}
\end{tikzpicture}
		\caption{The mean memory consumption depending on the number of variables, averaged over 10 systems.}
		\label{fig:mem}
	\end{subfigure}
\caption{The mean runtime and peak memory consumption for overdefined quadratic polynomial systems with increasing variable count $N$ and $M=2N$ polynomials over $\mathbb{F}_{101}$.}
\end{figure}
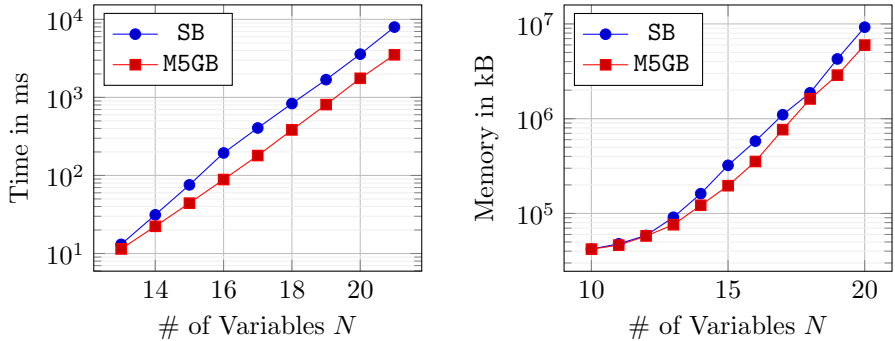

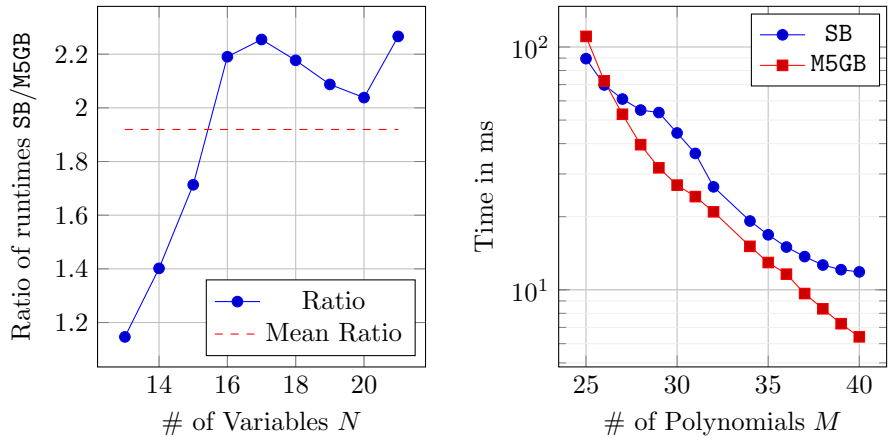
\begin{figure}
	\begin{subfigure}[t]{0.49\textwidth}
		\begin{tikzpicture}
    \pgfplotstabletranspose[input colnames to={n},colnames from={n}]{\mytabletoplot}{data/timing.csv}
    \begin{axis}[
        grid=both,
        major grid style = {lightgray},
        minor grid style = {lightgray!25},
        width=\textwidth,
        height=0.33\textheight,
        ylabel={Ratio of runtimes \texttt{SB}\slash\algo{}},
        xlabel={\# of Variables $N$},
        legend pos = south east,
        ytick distance = 0.2
        ]
        \addplot+[] table [x = {n}, y = {Ratio}] {\mytabletoplot};
        \addlegendentry{Ratio};
        \addplot+[no markers, dashed] table [x = {n}, y = {Mean}] 
        {\mytabletoplot};
        \addlegendentry{Mean Ratio};
    \end{axis}

\end{tikzpicture}
		  \caption{Ratio of runtimes between \texttt{SB} and \algo{}. The mean ratio is approximately 1.9.}
		  \label{fig:factor}
	\end{subfigure}
	\hfill
	\begin{subfigure}[t]{0.49\textwidth}
		\begin{tikzpicture}
    \pgfplotstabletranspose[input colnames to={n},colnames from={n}]{\mytabletoplot}{data/sanitized.csv}
    \begin{axis}[
        grid=both,
        major grid style = {lightgray},
        minor grid style = {lightgray!25},
        width=\textwidth,
        ymode = log,
        height=0.33\textheight,
        ylabel={Time in \SI{}{\milli\second}},
        xlabel={\# of Polynomials $M$},
        skip coords between index={8}{9},
        ]
        \addplot+[] table [x = {n}, y = {tm}] {\mytabletoplot};
        \addlegendentry{\texttt{SB}};
        \addplot+[] table [x = {n}, y = {tg}] {\mytabletoplot};
        \addlegendentry{\algo};

    \end{axis}
\end{tikzpicture}
		\caption{Runtime per basis element for an overdefined system over $\mathbb{F}_{101}$ with $N=15$ variables and a varying number of $M$ polynomials.}
		\label{fig:defined}
	\end{subfigure}
	\caption{Ratio of runtimes between \texttt{SB} and \algo{} for a varying variable count and runtimes of \texttt{SB} and \algo{} for a fixed number of variables $N$ and a varying number of polynomials $M$ over $\mathbb{F}_{101}$.}
\end{figure}

Figure \ref{fig:dq_timing1} illustrates the arithmetic mean of the measured timing results. The obtained results show that our implementation outperforms \texttt{SB} for dense quadratic systems in all tested systems. For both implementations, the time per basis element approximately doubles with each variable. Nevertheless, the runtime of \algo{} consistently stays below the runtime of \texttt{SB} in any of the tested systems. Our evaluation shows that the runtime ration between \texttt{SB} and \algo{} fluctuates over different values for $N$. Figure \ref{fig:factor} depicts the ratio between the arithmetic means of runtimes depending on the variable count $N$. After a slight drop between $N=17$ and $N=20$, for larger systems with $N=21$ the performance advantage of \algo{} starts to increase again.

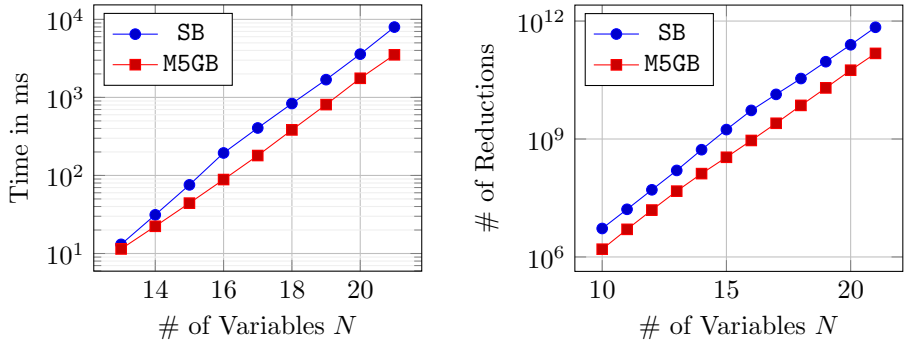
\begin{figure}[t]
	\begin{subfigure}[t]{0.49\textwidth}
		\begin{tikzpicture}
    \pgfplotstabletranspose[input colnames to={n},colnames from={n}]{\mytabletoplot}{data/timing.csv}
    \begin{axis}[
        grid=both,
        major grid style = {lightgray},
        minor grid style = {lightgray!25},
        width=\textwidth,
        ymode = log,
        ylabel={Time in \SI{}{\milli\second}},
        xlabel={\# of Variables $N$},
        legend pos= north west,
        ]
        \addplot+[] table [x = {n}, y = {SigGb}] {\mytabletoplot};
        \addlegendentry{\texttt{SB}};
        \addplot+[] table [x = {n}, y = {Chimera}] {\mytabletoplot};
        \addlegendentry{\algo};
    \end{axis}
\end{tikzpicture}
		\caption{
			The mean time spent per basis element calculated over 10 dense quadratic systems per variable count.
		}
		\label{fig:dq_timing2}
	\end{subfigure}
	\hfill
	\begin{subfigure}[t]{0.49\textwidth}
		\begin{tikzpicture}
    \pgfplotstabletranspose[input colnames to={n},colnames from={n}]{\mytabletoplot}{data/reductions.csv}
    \begin{axis}[
        grid=both,
        major grid style = {lightgray},
        minor grid style = {lightgray!25},
        width=\textwidth,
        ymode = log,
        legend pos= north west,
        ylabel={\# of Reductions},
        xlabel={\# of Variables $N$},
        ]
        \addplot+[] table [x = {n}, y = {SigGb}] {\mytabletoplot};
        \addlegendentry{\texttt{SB}};
        \addplot+[] table [x = {n}, y = {Chimera}] {\mytabletoplot};
        \addlegendentry{\algo};

    \end{axis}
\end{tikzpicture}
		\caption{The mean number of reductions steps for \texttt{SB} and \algo{}.}
		\label{fig:reductions}
	\end{subfigure}
\caption{Comparison between mean runtime spent per basis element and the mean number of reduction steps. The average is the arithmetic mean over 10 runs for dense, overdefined, quadratic polynomial systems over $\mathbb{F}_{101}$.}
\end{figure}

As we are particularly interested in dense quadratic systems, we also evaluate the performance depending on the number of polynomials $M$ for a fixed number of variables $N$. We find that decreasing the number of polynomials negatively influences the runtime and the performance gain of \algo{} compared to \texttt{SB}. For all evaluated systems, increasing the equation count reduces the runtime for both implementations. The runtime ratio is not strongly influenced by increasing the number of polynomials $M$. Figure \ref{fig:defined} illustrates the runtime changes depending on the number of provided polynomials $M$. Our baseline system has $N=15$ variables, and we vary the equation count. Our evaluation demonstrates that once $M \geq 1.8 \cdot N$ holds, our implementation outperforms the classic implementation for the tested systems. From these results, we conjecture that \algo{} will continue to outperform the classic algorithm for even larger systems, as long as they are sufficiently overdefined.    

Figure \ref{fig:mem} shows the memory consumption of both implementations on a logarithmic scale. As depicted, \algo{} tends to use less memory than the \texttt{SB} implementation. We can see that memory consumption increases exponentially. This is unsurprising, given that the number of possible S-pairs also grows exponentially. We could further reduce the memory footprint by design choices in the implementation. This, however, would not improve the memory growth behaviour of the algorithm itself. For the peak memory consumption metric, we were not able to test systems with $N=21$, as the memory profiling imposes an additional runtime overhead.  

In Figure \ref{fig:reductions} we see that the number of actually performed reductions is significantly lower in \algo{} than in \texttt{SB}. A comparison of \Cref{fig:dq_timing2} and \Cref{fig:reductions} indicates that the number of reductions is a good indicator of the time cost.
Note that the ratio between the reduction counts is larger than the actual ratio of runtimes in \Cref{fig:factor}. As our implementation is meant as a proof of concept, this observation leads to the assumption that a well-optimized implementation might lead to even higher performance gains.

\FloatBarrier

\section{Future Work}

As future work we plan to implement our \algo{} algorithm via the dedicated and optimized implementation of \texttt{M4GB}. It is the optimized implementation of \texttt{M4GB} that holds some of the top rankings in the MQ challenge \cite{mq}. Also, this would allow to draw a more stressable comparison between \algo{} and \texttt{M4GB}.

\bibliographystyle{alpha}
\bibliography{refs}

\end{document}